\documentclass[11pt, amsfonts]{amsart}

\usepackage{amsmath,amssymb,amsthm}
\usepackage{graphicx}
\usepackage[colorlinks=true]{hyperref}
\usepackage[all]{xy}

\numberwithin{equation}{section}

\SelectTips{eu}{12}

\newtheorem{theorem}{Theorem}[section]
\newtheorem{proposition}[theorem]{Proposition}
\newtheorem{lemma}[theorem]{Lemma}
\newtheorem{corollary}[theorem]{Corollary}

\theoremstyle{definition}
\newtheorem{definition}[theorem]{Definition}
\newtheorem{example}[theorem]{Example}

\theoremstyle{remark}
\newtheorem{remark}[theorem]{Remark}

\newcommand{ \Z}{\mathbb{Z}}
\newcommand{ \Q}{\mathbb{Q}}
\newcommand{ \C}{\mathbb{C}}

\newcommand{ \B}{\mathcal{B}}

\newcommand{ \U}{\operatorname{U}}

\newcommand{ \Gr}{\operatorname{Gr}}

\newcommand{ \rank}{\operatorname{rank}}

\newcommand{\ind}{\operatorname{ind}}
\newcommand{\prop}{\operatorname{prop}}
\newcommand{\trace}{\operatorname{trace}}

\title{Homotopy type of the space of finite propagation unitary operators on $\mathbb{Z}$}

\author{Tsuyoshi Kato}
\address{Department of Mathematics, Kyoto University, Kyoto, 606-8502, Japan}
\email{tkato@math.kyoto-u.ac.jp}

\author{Daisuke Kishimoto}
\address{Department of Mathematics, Kyoto University, Kyoto, 606-8502, Japan}
\email{kishi@math.kyoto-u.ac.jp}

\author{Mitsunobu Tsutaya}
\address{Faculty of Mathematics, Kyushu University, Fukuoka 819-0395, Japan}
\email{tsutaya@math.kyushu-u.ac.jp}



\thanks{Kato was supported by JSPS KAKENHI 17K18725
and 17H06461. Kishimoto was supported by JSPS KAKENHI 17K05248 and 19K03473. Tsutaya was supported by JSPS KAKENHI 19K14535}

\begin{document}

  \maketitle

  \baselineskip.525cm

  \begin{abstract}
    The index theory for the space of finite propagation unitary operators was developed by Gross, Nesme, Vogts and Werner from the viewpoint of quantum walks in mathematical physics. In particular, they proved that $\pi_0$ of the space is determined by the index. However, nothing is known about the higher homotopy groups. In this article, we describe the homotopy type of the space of finite propagation unitary operators on the Hilbert space of square summable $\C$-valued $\Z$-sequences, so we can determine its homotopy groups. We also study the space of (end-)periodic finite propagation unitary operators.
  \end{abstract}


  \section{Introduction}

  It is a remarkable fact that the wave operator on a complete Riemannian manifold has the finite propagation speed. More precisely, for a complete Riemannian manifold $X$ and the Dirac operator $D$ acting on spinor bundles $S^{\pm}$, the wave operator $e^{\sqrt{-1}tD}$ acts on $L^2(X; S^{\pm})$ such that its kernel function $k_t(x,y)$ has the support within $|t|$ neighborhood of the diagonal. Then we can get the asymptotic expansion of the kernel function $k_t(x,y)$ by using a finite propagation unitary operator, and so finite propagation unitary operators are basic tools for analyzing important objects in mathematical physics. See \cite{R} for details. Several other mathematical and physical phenomena with finite propagation speed can be analyzed by using finite propagation unitary operators as well. Then finite propagation unitary operators have been studied intensely in several contexts \cite{C,CGT,GNVW,R}. However little is known about the topology of finite propagation unitary operators.

  Gross, Nesme, Vogts and Werner \cite{GNVW} studied the topology of finite propagation unitary operators in the context of quantum walks in mathematical physics, where their quantum walks are presented by finite propagation unitary operators. They introduced the integer valued index of finite propagation unitary operators and investigated its properties. In particular, they proved that the index classifies the connected components of the space of finite propagation unitary operators, that is, its $\pi_0$ is isomorphic with $\Z$ through the index. This indicates that the space of finite propagation unitary operators has non-trivial homotopy type while the space of all unitary operators is contractible \cite{Ku}. However, nothing is known about its higher homotopy groups and its homotopy type. Thus it is of particular importance to study the homotopy groups and the homotopy type of the space of finite propagation unitary operators.

  In this article, we study the homotopy groups and the homotopy type of the space of finite propagation unitary operators $\mathcal{U}$ on the Hilbert space $\ell^2(\C)$ of square summable $\C$-valued $\Z$-sequences. The topology of $\mathcal{U}$ is defined as follows. Let $\mathcal{U}_L$ denote the set of unitary operators on $\ell^2(\C)$ of propagation $L$. We equip $\mathcal{U}_L$ with the norm topology. We have $\mathcal{U}=\bigcup_{L>0}\mathcal{U}_L$ as a set, and we equip $\mathcal{U}$ with the weak topology with respect to $\mathcal{U}_L$. First, we will compute the homotopy groups of $\mathcal{U}$.

  \begin{theorem}
    \label{main 1}
    There is an isomorphism
    \[
      \pi_i(\mathcal{U})\cong
      \begin{cases}
        \ell^\infty(\Z)_S & \text{$i$ is odd}\\
        \mathbb{Z} & \text{$i$ is even.}
      \end{cases}
    \]
  \end{theorem}

  Remarks on the abelian group $\ell^\infty(\Z)_S$ in Theorem \ref{main 1} are in order. Let $\ell^\infty(\Z)$ denote the abelian group of bounded $\Z$-valued $\Z$-sequences, and let $S\colon\ell^\infty(\Z)\to\ell^\infty(\Z)$ denote the left shift operator defined by $S(v_i)_i =(v_{i+1})_i$. The abelian group $\ell^\infty(\Z)_S$ denotes the group of coinvariants of the action of $S$ on $\ell^\infty(\Z)$, that is,
  \[
    \ell^\infty(\Z)_S=\ell^\infty(\Z)/\{ a-Sa\mid a\in\ell^\infty(\Z)\}.
  \]
  We will prove in Proposition \ref{coinvariant} that $\ell^\infty(\Z)_S$ is a $\mathbb{Q}$-vector space, where the dimension of $\ell^\infty(\Z)_S$ is uncountable because it is an uncountable set.

  Next, we will describe the homotopy type of $\mathcal{U}$. In the computation of the homotopy groups, we will consider the homotopy fibration involving $\mathcal{U}$. We will prove the homotopy fibration is trivial by using the above result that the odd homotopy groups of $\mathcal{U}$ are $\mathbb{Q}$-vector spaces together with a result on the homotopy type of the space of periodic finite propagation unitary operators that we will obtain previously. This will give us a description of the homotopy type of $\mathcal{U}$ in terms of familiar spaces. To state the result, we set notation. Let
  \[
  	\U(\infty)=\lim_{n\to\infty}\U(n)
  \]
  and let $BG$ denote the classifying space of a topological group $G$. Let $K(A,n)$ denote the Eilenberg--MacLane space of type $(A,n)$ for an abelian group $A$ and a positive integer $n$. For pointed spaces $\{X_n\}_{n\ge 1}$, let
  \[
  	\overset{\circ}{\prod}_{n\ge 1} X_n=\lim_{N\to\infty}\prod_{n=1}^N X_n
  \]
  which is called the weak product of $X_n$. Now we state the main theorem of this article.

  \begin{theorem}
  \label{main 2}
  	There is a homotopy equivalence
  	\[
      \mathcal{U}\simeq\mathbb{Z}\times B\U(\infty)\times\overset{\circ}{\prod_{n\ge 1}} K(\ell^\infty(\Z)_S,2n-1).
    \]
  \end{theorem}

  As mentioned above, we will also consider the subgroups of $\mathcal{U}$ consisting of (end-)periodic finite propagation unitary operators. We will determine their homotopy types, which will be used to prove the triviality of the homotopy fibration involving $\mathcal{U}$.

  There are at least two directions of further research. One is to construct a bundle theory and the obstruction theory for bundles whose structure groups are the group of finite propagation unitary operators. The other is to generalize our research to the completion of the finite propagation unitary operators, which is the unitary group of an example of a Roe algebra. The former topic will be developed as our next subject, and the latter is studied in \cite{KKT}.

  \subsection*{Acknowledgement}

  The authors are grateful to the anonymous referee and the corresponding editor Jonathan Rosenberg for useful advice and comments.


  \section{Index theory}
  \label{section_index}

  In this section, we recall from \cite{GNVW} the index theory of finite propagation unitary operators and its related properties that we are going to use.


  \subsection{Finite propagation unitary operators}

  Let $\ell^2(\mathbb{C})$ denote the Hilbert space of square summable $\C$-valued $\Z$-sequences, and let $\B$ denote the space of bounded operators on $\ell^2(\C)$ equipped with the norm topology. We define finite propagation unitary operators on $\ell^2(\mathbb{C})$. Any bounded operator $T$ on $\ell^2(\mathbb{C})$ can be expressed by the matrix form as
  \[
    T=(T_{ij})
  \]
  where $i$ and $j$ range over all integers. Then we will discuss bounded operators on $\ell^2(\mathbb{C})$ in terms of the matrix form.

  \begin{definition}
    We say that a bounded operator $T=(T_{ij})\in \B$ has \textit{finite propagation} if the propagation
    \[
      \prop(T)=\sup\{ |i-j|\mid T_{ij}\ne 0 \}
    \]
    is finite.
  \end{definition}

  By an obvious property of matrix multiplication, for any finite propagation operators $S,T\in\mathcal{B}$, the following inequality holds:
    \[
      \prop(ST) \le  \prop(S)+\prop(T).
    \]
  Then in particular, the set of finite propagation operators is closed under product.

  \begin{definition}
    Let $\mathcal{U}_L$ denote the space of unitary operators of propagation $\le L$, which is equipped with the norm topology. Define the space $\mathcal{U}$ as the direct limit
    \[
      \mathcal{U}=\lim_{L\to\infty}\mathcal{U}_L
    \]
    which is equipped with the weak topology with respect to $\mathcal{U}_L$.
  \end{definition}

  While $\mathcal{U}_L$ is not closed under product, $\mathcal{U}$ is closed under product, so it is a topological group.

  \begin{remark}
    Though the inclusion $\mathcal{U}_L\to \B$ induces the injective continuous map $\mathcal{U}\to \B$, it is not a homeomorphism onto its image. Namely, the topology of $\mathcal{U}$ that we consider here is different from the subspace topology of $\B$.
  \end{remark}


  \subsection{Index}

  We define the index of finite propagation unitary operators. We have another matrix expression for $T\in \B$ by
  \[
    T=
    \begin{pmatrix}
      T_{--} & T_{-+} \\
      T_{+-} & T_{++}
    \end{pmatrix},
  \]
  where for $T=(T_{ij})$
  \[
    T_{--}=(T_{ij})_{i<0,j<0},\quad T_{-+}=(T_{ij})_{i<0,j\ge 0},\quad T_{+-}=(T_{ij})_{i\ge 0,j<0},\quad T_{++}=(T_{ij})_{i\ge 0,j\ge 0}.
  \]

  \begin{definition}
    The \emph{index}  of $U\in\mathcal{U}$ is defined by
    \[
      \ind(U)=\|U_{-+}\|_{\mathrm{HS}}^2-\|U_{+-}\|_{\mathrm{HS}}^2
    \]
    where $\|\cdot\|_{\mathrm{HS}}$ denotes the Hilbert--Schmidt norm, that is, $\|T\|_{\mathrm{HS}}^2=\sum_{ij}|T_{ij}|^2$ for $T=(T_{ij})$.
  \end{definition}

  The index of a finitely generated unitary operator $U$ on $\ell^2(\mathbb{C})$ is well-defined because the components $U_{-+}$ and $U_{+-}$ have only finitely many non-zero entries.

  \begin{example}
    The \emph{left shift operator} $S\in\mathcal{B}$ is defined by
    \[
      S(v_i)_{i}=(v_{i+1})_{i}.
    \]
    The operator $S$ is a unitary operator with $\prop(S)=1$ and $\ind(S)=1$. This finite propagation unitary operator will play an important role later in our study.
  \end{example}

  We recall from \cite{GNVW} properties of the index of finite propagation unitary operators that we are going to use later on.

  \begin{theorem}
    \cite[Theorem 3]{GNVW}
    \label{index}
    \begin{enumerate}
      \item For any $U\in\mathcal{U}$, $\ind(U)$ is an integer.

      \item For $U,V\in\mathcal{U}$, the additivity holds:
      \[
        \ind(UV)=\ind(U)+\ind(V).
      \]

      \item The index defines a homomorphism $\ind\colon\mathcal{U}\to\mathbb{Z}$, which induces an isomorphism
      \[
        \pi_0(\mathcal{U})\cong \mathbb{Z}.
      \]
    \end{enumerate}
  \end{theorem}

  A slight generalization of the item 4 of \cite[Theorem 3]{GNVW} will be used to construct a certain homogeneous space related to $\mathcal{U}$. The proof of it in \cite{GNVW} is rather implicit, so we here give a direct proof. Let $\{e_i\}_{i\in\mathbb{Z}}$ denote the standard orthonormal basis of $\ell^2(\mathbb{C})$ and $P\colon\ell^2(\mathbb{C})\to\ell^2(\mathbb{C})$ be the projection defined by
  \[
    Pe_i=
    \begin{cases}
      e_i & i\ge 0\\
      0 & i<0.
    \end{cases}
  \]
  Let $Q\colon\ell^2(\mathbb{C})\to\ell^2(\mathbb{C})$ be any orthogonal projection satisfying
  \[
    Qe_i  \ = \
    \begin{cases}
      e_i & i\ge k+L\\
      0 & i < k
    \end{cases}
  \]
  for some integer $k$ and some positive integer $L$. Then we can define the trace of $P-Q$ because $P-Q$ has only finitely many non-zero entries.

  \begin{lemma}
    \label{decomp lem}
    If $\trace(P-Q)=0$, then there exists a unitary operator $V$ on $\ell^2(\mathbb{C})$ such that
  	\[
  		P=V^\ast QV\quad\text{and}\quad Ve_i=e_i\quad\text{unless}\quad k\le i<k+L.
  	\]
  \end{lemma}

  \begin{proof}
  Let $W$ denote the subspace of $\ell^2(\mathbb{C})$ spanned by $e_k,e_{k+1},\ldots, e_{k+L-1}$. Then we have
  \[
  	P(W)\subset W\quad\text{and}\quad Q(W)\subset W
  \]
  or equivalently
  \[
  	P(W^{\perp})\subset W^{\perp}\quad\text{and}\quad Q(W^{\perp})\subset W^{\perp}.
  \]
  It follows from $\trace(P-Q)=0$ that $k+L \ge 0$ and $k \le 0$, and so $P=Q$ on $W^{\perp}$. Then
  \[
  0=\trace(P-Q)=\trace(P-Q)|_W=\trace(P|_W-Q|_W)=\trace(P|_W)-\trace(Q|_W)
  \]
  implying $\rank(Q|_W)=\rank(P|_W)$, where $W$ is a finite dimensional vector space so that we can consider traces and ranks of linear maps on $W$.

  Now we take an orthonormal basis $v_0,v_1,\ldots,v_{k+L-1}$ of $Q(W)$ and extend it to an orthonormal basis $v_k,v_{k+1},\ldots,v_{k+L-1}$ of $W$, where $k+L \ge 0$ and $k \le 0$ as above. By definition,
  \[
    Qv_i  \ = \
    \begin{cases}
      v_i & i\ge 0\\
      0 & i < 0.
    \end{cases}
  \]
  Define the unitary operator $V$ on $\ell^2(\mathbb{Z})$ by
  \[
  	Ve_i=
  	\begin{cases}
  		v_i&k\le i<k+L\\
  		e_i&\text{otherwise}.
  	\end{cases}
  \]
  Then $P=V^*QV$, so the proof is complete.
  \end{proof}

  We consider the $C^\ast$-algebra of block diagonal operators
  \[
  B_k(L)=\{ T\in \B\mid T_{ij}=0\text{ unless $k+nL\le i,j<k+(n+1)L$ for some $n\in\mathbb{Z}$}\}.
  \]
  for $k\in\mathbb{Z}$, $L\ge1$, which actually depends only on $L$ and $k$ mod $L$. By definition, an element of $B_k(L)$ has the matrix form illustrated below.
  \begin{center}
    \includegraphics[scale=0.8]{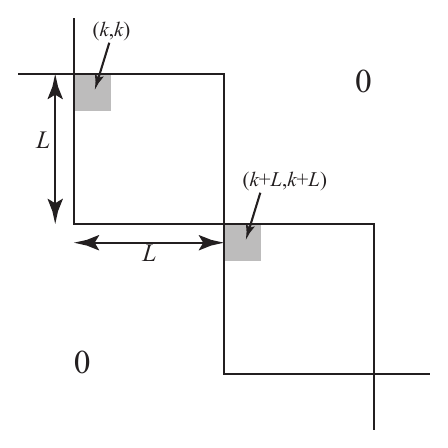}
  \end{center}
  We will write an element of $B_k(L)$ as $U=(U_j)_{j\in\mathbb{Z}}$, where $U_j$ is the block containing the $(k+jL,k+jL)$-entry.
  For a unital $\ast$-algebra $A$, let
  \[
  	\U(A)=\{x\in A\mid xx^\ast=x^\ast x=1\}.
  \]
  Then we can consider $\U(B_k(L))$, for instance. Now we give a direct proof of a slight extension of the item 4 of \cite[Theorem 3]{GNVW}.

  \begin{theorem}
  	\label{decomp}
  	If $U\in\mathcal{U}$ satisfies $\ind(U)=0$ and $\prop(U)\le L$, then it admits a product decomposition
  	\[
  		U=VW
  	\]
  	for some $V \in\U(B_0(2L))$ and $W \in\U(B_{-L}(2L))$.
  \end{theorem}

  \begin{proof}
  Clearly, $U^\ast PU$ is an orthogonal projection such that
  \[
  	U^\ast PUe_i=
  	\begin{cases}
  		e_i&i\ge L\\
  		0&i<-L.
  	\end{cases}
  \]
  Then it follows that
  \[
    P-U^\ast PU=
    \begin{pmatrix}
      0 & 0 \\
      0 & I
    \end{pmatrix}
    -
    \begin{pmatrix}
      U_{--}^\ast & U_{+-}^\ast \\
      U_{-+}^\ast & U_{++}^\ast
    \end{pmatrix}
    \begin{pmatrix}
      0 & 0 \\
      0 & I
    \end{pmatrix}
    \begin{pmatrix}
      U_{--} & U_{-+} \\
      U_{+-} & U_{++}
    \end{pmatrix}
    =
    \begin{pmatrix}
      -U_{+-}^\ast U_{+-} & -U_{+-}^\ast U_{++} \\
      -U_{++}^\ast U_{+-} & U_{-+}^\ast U_{-+}
    \end{pmatrix}
  \]
  implying
  \[
   \trace(P-U^\ast PU)=\|U_{-+}\|_{\mathrm{HS}}^2-\|U_{+-}\|_{\mathrm{HS}}^2=\ind(U)=0.
  \]
  Thus we can apply Lemma \ref{decomp lem} to get a unitary operator $\widetilde{V}_0$ on $\ell^2(\mathbb{C})$ such that
  \[
  	P=\widetilde{V}_0^\ast U^\ast PU\widetilde{V}_0\quad\text{and}\quad\widetilde{V}_0e_i=e_i\quad\text{unless}\quad -L\le i<L.
  \]
  These are interpreted in the matrix form as follows. Since $U\widetilde{V}_0$ commutes with the projection $P$, the $(-,+)$ and $(+,-)$ components of $U\widetilde{V}_0$ are trivial. Since $\widetilde{V}_0e_i=e_i$ unless $-L\le i<L$, the non-trivial block component of $\widetilde{V}_0$
  is a $2L\times 2L$ matrix.

  Applying the same argument to $S^{2nL}US^{-2nL}$ instead of $U$, which is the diagonal $2n$-shift of $U$, we obtain $2L\times 2L$ matrices $W_n$ for $n\in\Z$ such that for
  \[
    W=
    \begin{pmatrix}
      \ddots \\
      & W_{-1}^\ast \\
      & & W_0^\ast \\
      & & & W_1^\ast \\
      & & & & \ddots
    \end{pmatrix}
    \in\U(B_{-L}(2L))
  \]
  we have $UW^\ast\in\U(B_0(2L))$, completing the proof.
  \end{proof}
  \begin{remark}
  \label{rem_decomp_periodic}
  As in the proof, the operator $V$ in Lemma \ref{decomp lem} is actually chosen depending only on the information of $Q$ around $(0,0)$.
  This implies that when $U$ is periodic as we consider in Section \ref{section_periodic}, the matrices $W_n$ ($n\in\mathbb{Z}$) in the proof of the theorem can be chosen as the same ones.
  As a consequence, $V$ is also a diagonal block operator with the same diagonal blocks. 
  \end{remark}


  \section{Homotopy groups of $\mathcal{U}$}

  In this section, we compute the homotopy groups of $\mathcal{U}$. The main ingredient in the computation is the space $\mathcal{W}_L$ which approximates the space $\mathcal{U}$. We will construct a principal bundle involving the approximating space $\mathcal{W}_L$, and will analyze it to compute the homotopy groups of $\mathcal{U}$.


  \subsection{Homotopy groups of $\U(B_k(L))$}

  We compute the homotopy groups of the space of diagonal block unitary operators $\U(B_k(L))$, which will be needed to compute the homotopy groups of $\mathcal{U}$.

  Let $\U(L)$ denote the rank $L$ unitary group, and let $X^\Z$ denote the product of countably infinite copies of a space $X$. It is clear that $\U(B_k(L))$ and $\U(L)^\Z$ are the same as sets, where the homotopy groups of $\U(L)$ are well known in dimension $\le 2L-1$. However, their topologies are different so that the homotopy groups of $\U(B_k(L))$ are not directly computed from those of $\U(L)^\Z$. Then we need a closer look at the topology of $\U(B_k(L))$ to compute its homotopy groups.

  Let
  \[
    p_j\colon \U(B_k(L))\to\U(L)
  \]
  denote the projection onto the block component containing the $(k+jL,k+jL)$-entry. We consider a control on maps from spheres into $B_k(L)$ by the maps $p_j$. Recall that a family of maps $\{f_i: X \to Y\}_{i \in \mathbb {Z}} $ between  metric spaces is called uniformly equicontinuous if for any $\epsilon >0$, there is $\delta >0$ such that $d_Y(f_n(x),f_n(x')) < \epsilon$ holds for any $n \in \mathbb {Z}$ and $x,x'\in X$ with $d_X(x,x') < \delta$. We say that a map, not necessarily continuous, $f\colon X\to\U(B_k(L))$ is of finite type if there are finitely many continuous maps $g_1,\ldots,g_k\colon X\to\U(L)$ such that for each $j\in\Z$, $p_j\circ f\colon X\to\U(L)$ is one of $g_1,\ldots,g_k$. For example, given a continuous map $X\to\U(L)$, the diagonal map
  \[
    X\to\U(B_k(L)),\quad x\mapsto(\ldots,f(x),f(x),f(x),\ldots)
  \]
  is of finite type.

  \begin{lemma}
    \label{ascoli-arzela}
    \begin{enumerate}
	\item
      Any map $f\colon S^n\to\U(B_k(L))$ of finite type is continuous.
	\item
      Any map $f\colon S^n\to\U(B_k(L))$ is homotopic to a map of finite type.
    \end{enumerate}
  \end{lemma}

  \begin{proof}
    (1) Let $f\colon S^n\to\U(B_k(L))$ be of finite type. Then there are finitely many continuous map $g_1,\ldots,g_k\colon S^n\to\U(L)$ such that for each $j\in\Z$, $p_j\circ f$ is one of $g_1,\ldots,g_k$. So we get
    \[
      \|f(x)-f(y)\|\le\sup\{\|p_j\circ f(x)-p_j\circ f(y)\|\mid j\in\Z\}=\max\{\|g_i(x)-g_i(y)\|\mid 1\le i\le k\}.
    \]
    Thus since each $g_i$ is continuous, $f$ is continuous.

    (2) Since $\|T\|=\sup_j\|p_j\circ T\|$ for any $T\in \U(B_k(L))$ and $S^n$ is compact, the family of maps $\{p_j\circ f\}_j$ is uniformly equicontinuous. Then since $\U(L)$ is compact bounded in $M_L(\mathbb{C})\cong\mathbb{C}^{L^2}$, it follows from the Arzel\`{a}--Ascoli theorem that the family of maps $\{p_j\circ f\}_j$ is relatively compact in the space of maps $S^n\to\U(L)$. Hence there is a finite collection of maps $h_1,\ldots,h_r\colon S^n\to\U(L)$ such that any $p_j\circ f$ is sufficiently close to one of $h_1,\ldots,h_r$. Thus there are homotopies between these maps by using geodesics in $\U(L)$ or the local path-connectedness of the space of unitaries $\U(C(S^n)\otimes B_k(L))$ of the tensor product $C(S^n)\otimes B_k(L)$ as $C^\ast$-algebra, which yields a homotopy from $f$ to the desired map $S^n\to\U(B_k(L))$.
  \end{proof}

  Now we compute the homotopy groups of $\U(B_k(L))$. Let $\ell^\infty(\mathbb{Z})$ denote the group of bounded $\mathbb{Z}$-valued $\mathbb{Z}$-sequences.

  \begin{proposition}
  	\label{pi(B)}
  	For $i\le 2L-1$, there is an isomorphism
  	\[
  		\pi_i(\U(B_k(L)))\cong
  		\begin{cases}
  			\ell^\infty(\Z)& i\text{ is odd}\\
  			0 & i\text{ is even}.
  		\end{cases}
  	\]
  \end{proposition}

  \begin{proof}
  	By Lemma \ref{ascoli-arzela} (2), a map $f\colon S^n\to\U(B_k(L))$ is null-homotopic if and only if $p_j\circ f$ is null-homotopic for each $j$. Then the map $\rho=\prod_{j\in\Z}p_j\colon\U(B_k(L))\to\U(L)^{\mathbb{Z}}$ is injective in homotopy groups. It remains to determine the image of the map $\rho$ in homotopy groups of dimension $\le 2L-1$. For $i\le 2L-1$, there is an isomorphism
  	\[
  	\pi_i(\U(L))\cong
  	\begin{cases}
  		\Z& i\text{ is odd}\\
  		0 & i\text{ is even}.
  	\end{cases}
  	\]
  	Note that there is a natural isomorphism $\pi_*(\U(L)^\Z)\cong\pi_*(\U(L))^\Z$. Then $\pi_{2i}(\U(L)^\Z)=0$ for $i\le L$, implying $\rho_*(\pi_{2i}(\U(B_k(L)))=0$ for $i\le L$. On the other hand, $\pi_{2i-1}(\U(L)^\Z)$ is identified with the set of all $\Z$-valued $\Z$-sequences. Observe that an integer sequence $(a_i)_{i\in\Z}$ belongs to $\ell^\infty(\mathbb{Z})$ if and only if there are finitely many integers $b_1,\ldots,b_k$ such that for each $i\in\Z$ $a_i$ is one of $b_1,\ldots,b_k$.  Thus by Lemma \ref{ascoli-arzela} (1), $\rho_*(\pi_{2i-1}(\U(B_k(L)))$ for $i\le L$ is exactly $\ell^\infty(\Z)$, completing the proof.
  \end{proof}


  \subsection{The space $\mathcal{W}_L$}

  We define the space $\mathcal{W}$ which approximates $\mathcal{U}$. Let the group $\U(B_0(L))$ act on $\U(B_0(2L))\times\U(B_{-L}(2L))$ from the right by
  \[
    (U_1,U_2)\cdot V=(U_1V,V^{-1}U_2)
  \]
  for $(U_1,U_2)\in\U(B_0(2L))\times\U(B_{-L}(2L))$ and $V\in\U(B_0(L))$. We define the space $\mathcal{W}_L$ as the orbit space of this action. There is a natural inclusion $\mathcal{W}_L\to\mathcal{W}_{3L}$ as depicted below.
  \begin{center}
    \includegraphics[scale=0.8]{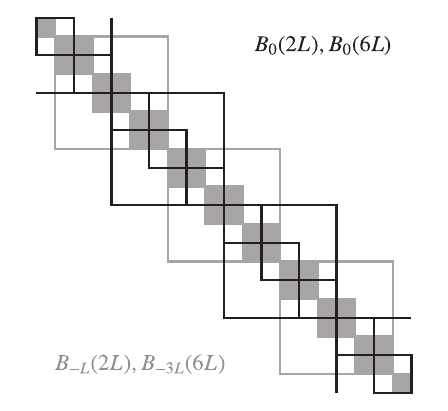}
  \end{center}
  We define the space
  \[
    \mathcal{W}=\lim_{n\to\infty}\mathcal{W}_{3^n}
  \]
  which is given the weak topology with respect to $\mathcal{W}_{3^n}$. 

  The following principal bundle is fundamental in investigating the space $\mathcal{W}_L$.

  \begin{lemma}
    \label{fibration W}
    For $L\ge1$, there is a principal bundle
    \[
      \U(B_0(L))\to \U(B_0(2L))\times\U(B_{-L}(2L))\to\mathcal{W}_L.
    \]
  \end{lemma}

  \begin{proof}
    By definition, $\U(B_0(2L))\times\U(B_{-L}(2L))$ is a $\U(B_0(L))$-torsor over $\mathcal{W}_L$. Then it remains to show the local triviality. The spaces $\U(B_0(L))$ and $\U(B_0(2L))\times\U(B_{-L}(2L))$ are Banach manifolds such that the map $\U(B_0(L))\to \U(B_0(2L))\times\U(B_{-L}(2L))$ is a closed embedding onto a direct summand at each tangent space. Then $\mathcal{W}_L$ has a natural Banach manifold structure such that the sequence in the statement is a principal bundle.
  \end{proof}

  Let $\mathcal{U}^0$ denote the identity component of $\mathcal{U}$. By Theorem \ref{index}, $\mathcal{U}^0$ is the subspace of $\mathcal{U}$ consisting of index zero elements.

  \begin{proposition}
    \label{U^0}
    There is a homeomorphism
    \[
      \mathcal{W}\cong\mathcal{U}^0.
    \]
  \end{proposition}

  \begin{proof}
    By Theorem \ref{decomp}, each $U\in\mathcal{U}^0$ with $\prop(U)\le L$ admits a product decomposition
    \[
      U=U_1U_2
    \]
    for $U_1\in\U(B_0(2L))$ and $U_2\in\U(B_{-L}(2L))$. Then the map
    \[
      \Phi_L\colon\mathcal{W}_L\to\mathcal{U}^0,\quad[U_1,U_2]\mapsto U_1U_2
    \]
    is continuous and injective. Since $\Phi_L$ is natural with respect to $L$, we get a map
    $\Phi\colon\mathcal{W}\to\mathcal{U}_0$. Since $\mathcal{W}$ is equipped with the weak topology, $\Phi$ is continuous, and by Theorem \ref{decomp}, $\Phi$ is bijective too. Then it remains to show that $\Phi$ is a closed map.

    By definition, the principal bundle of Lemma \ref{fibration W} is smooth such that the image of the differential of the fiber inclusion at each point is a closed direct summand. Then by the implicit function theorem, the space $\mathcal{W}_L$ is a smooth Banach manifold such that the principal bundle of Lemma \ref{fibration W} is smooth. In particular, the quotient map $\U(B_0(2L))\times\U(B_{-L}(2L))\to\mathcal{W}_L$ has a local section which is a closed map. On the other hand, the map
    \[
      \U(B_0(2L))\times\U(B_{-L}(2L))\to\mathcal{U}^0,\quad(U_1,U_2)\mapsto U_1U_2
    \]
    is a closed map. Thus the map $\Phi_L\colon\mathcal{W}_L\to\mathcal{U}^0$ turns out to be a closed map. Since $\mathcal{W}$ is given the weak topology, the map $\Phi$ is a closed map too, completing the proof.
  \end{proof}

  \begin{corollary}
    \label{pi(U^0)}
    There is an isomorphism
    \[
      \pi_*(\mathcal{U}^0)=\lim_{n\to\infty}\pi_*(\mathcal{W}_{3^n}).
    \]
  \end{corollary}

  \begin{proof}
    Since the inclusion $\mathcal{W}_L\to\mathcal{W}_{3L}$ is a cofibration for each $L$, we have
    \[
      \pi_*(\mathcal{W})=\lim_{n\to\infty}\pi_*(\mathcal{W}_{3^n})
    \]
    Then the proof is finished by Proposition \ref{U^0}.
  \end{proof}

  Thus we compute the homotopy groups of $\mathcal{W}_L$. To this end, we need the following lemma.

  \begin{lemma}
    \label{inclusion B}
    The inclusion $\U(B_0(L))\to\U(B_0(nL))$ is identified with the map
    \[
      \ell^\infty(\Z)\to\ell^\infty(\Z),\quad(a_i)_i\mapsto(a_{ni}+a_{ni+1}+\cdots+a_{ni+n-1})_i
    \]
    in $\pi_{2i-1}$ for $i\le L$.
  \end{lemma}

  \begin{proof}
    In the proof of Proposition \ref{pi(B)}, it is shown that the canonical map $\U(B_0(kL))\to\U(kL)^\Z$ is identified with the natural inclusion $\ell^\infty(\Z)\to\Z^\Z$ in $\pi_{2i-1}$ for $i\le L$. There is a commutative diagram
    \[
      \xymatrix{
        \U(B_0(L))\ar[r]\ar[d]&\U(L)^\Z\ar[d]\\
        \U(B_0(nL))\ar[r]&\U(nL)^\Z
      }
    \]
    where the right vertical map is the product of the diagonal block inclusion. Since the inclusion $\U(L)^\Z\to\U(nL)^\Z$ is identified with
    \[
      \Z^\Z\to\Z^\Z,\quad\quad(a_i)_i\mapsto(a_{ni}+a_{ni+1}+\cdots+a_{ni+n-1})_i
    \]
    in $\pi_{2i-1}$ for $i\le L$, the proof is complete.
  \end{proof}

  We compute the homotopy groups of $\mathcal{W}_L$.

  \begin{proposition}
    \label{pi(W)}
      For $1\le i\le 2L$, there is an isomorphism
      \[
        \pi_i(\mathcal{W}_L)=
        \begin{cases}
          \ell^\infty(\mathbb{Z})_S & \text{$i$ is odd}\\
          \mathbb{Z} & \text{$i$ is even.}
        \end{cases}
      \]
  \end{proposition}

  \begin{proof}
    Consider the homotopy exact sequence
    \[
      \cdots\to\pi_*(\U(B_0(L)))\to\pi_*(\U(B_0(2L)))\times\pi_*(\U(B_{-L}(2L)))\to\pi_*(\mathcal{W}_L)\to\cdots
    \]
    associated to the principal bundle in Lemma \ref{fibration W}. By Proposition \ref{pi(B)}, we get the exact sequence
    \[
      0 \to\pi_{2i}(\mathcal{W}_L) \to\ell^\infty(\mathbb{Z}) \xrightarrow{\alpha}\ell^\infty(\mathbb{Z})\times\ell^\infty(\mathbb{Z}) \to\pi_{2i-1}(\mathcal{W}_L) \to 0
    \]
    for $1\le i\le L$. By the construction, the fiber inclusion of the principal bundle is given by
    \[
      \U(B_0(L))\to\U(B_0(2L))\times\U(B_{-L}(2L)),\quad U\mapsto(U,U^{-1}).
    \]
    Then by Lemma \ref{inclusion B}, the map $\alpha$ is given by
    \begin{equation}
      \label{alpha def}
      \alpha((a_j)_j)=((a_{2j}+a_{2j+1})_j,(-a_{2j-1}-a_{2j})_j)
    \end{equation}
    for $(a_j)_j\in\ell^\infty(\Z)$, where the inclusion $\U(B_0(L))\to\U(B_{-L}(2L))$ is identified with the inclusion $\U(B_0(L))\to\U(B_{0}(2L))$ by the diagonal shift. Thus
    \[
      \pi_{2i}(\mathcal{W}_L)\cong\mathrm{Ker}\,\alpha=\{(\ldots,a,-a,a,-a,\ldots)\mid a\in\Z\}\cong\Z.
    \]
    Consider the isomorphism
    \[
      \beta\colon\ell^\infty(\Z)\times\ell^\infty(\Z)\to\ell^\infty(\Z),\quad((a_j)_j,(b_j)_j)\mapsto(a_j,b_j)_j.
    \]
    Then the image of the map $\beta\circ\alpha\colon\ell^\infty(\Z)\to\ell^\infty(\Z)$ consists of elements of the form
    \begin{align*}
      &(\ldots,a_{-2}+a_{-1},-a_{-3}-a_{-2},a_0+a_1,-a_{-1}-a_0,a_2+a_3,-a_1-a_2,\ldots)\\
      &=(\ldots,a_{-2},-a_{-2},a_0,-a_0,a_2,-a_2,\ldots)+(\ldots,a_{-1},-a_{-3},a_1,-a_{-1},a_3,-a_1,\ldots)\\
      &=-(1-S)(\ldots,0,a_{-2},0,a_0,0,a_2,\ldots)-(1-S^3)(\ldots,0,a_{-1},0,a_1,0,a_3,0,\ldots)\\
      &=-(1-S)((\ldots,0,a_{-2},0,a_0,0,a_2,\ldots)+(1+S+S^2)(\ldots,0,a_{-1},0,a_1,0,a_3,0,\ldots))
    \end{align*}
    where $S$ denotes the shift operator. In the last line, $a_0$ in the first term and $a_1$ in the second term are placed at the zeroth and the first entry, respectively. It is easy to see that any element of $\ell^\infty(\Z)$ can be expressed by an element of the form
    \[
      (\ldots,b_{-2},0,b_0,0,b_2,0,\ldots)+(1+S+S^2)(\ldots,b_{-1},0,b_1,0,b_3,0,\ldots).
    \]
    Thus $\mathrm{Im}\,\beta\circ\alpha=\mathrm{Im}(1-S)$, and therefore $\pi_{2i}(\mathcal{W}_L)\cong\ell^\infty(\Z)_S$ for $1\le i\le L$, completing the proof.
  \end{proof}

  Now we are ready to prove Theorem \ref{main 1}.

  \begin{proof}
    [Proof of Theorem \ref{main 1}]
    The inclusion $\mathcal{W}_L\to\mathcal{W}_{3L}$ is natural with respect to the principal bundle of Lemma \ref{fibration W} in the sense that there is a commutative diagram
    \[
      \xymatrix{
        \U(B_0(L)) \ar[r] \ar[d]& \U(B_0(2L))\times\U(B_{-L}(2L)) \ar[r] \ar[d]& \mathcal{W}_L \ar[d] \\
        \U(B_0(3L)) \ar[r]& \U(B_0(6L))\times\U(B_{-3L}(6L)) \ar[r]& \mathcal{W}_{3L}.
      }
    \]
    Then by Proposition \ref{pi(B)}, we get a commutative diagram for $1\le i\le L$
    \[
    \xymatrix{
    0 \ar[r]& \pi_{2i}(\mathcal{W}_L) \ar[r]^{\theta_0} \ar[d]& \ell^\infty(\mathbb{Z}) \ar[r]^(.4){\alpha_0} \ar[d]^{\phi_1}& \ell^\infty(\mathbb{Z})\times\ell^\infty(\mathbb{Z}) \ar[r] \ar[d]^{\phi_2}& \pi_{2i-1}(\mathcal{W}_L) \ar[r] \ar[d]& 0\\
    0 \ar[r]& \pi_{2i}(\mathcal{W}_{3L}) \ar[r]^{\theta_1}& \ell^\infty(\mathbb{Z}) \ar[r]^(.4){\alpha_1}& \ell^\infty(\mathbb{Z})\times\ell^\infty(\mathbb{Z}) \ar[r]& \pi_{2i-1}(\mathcal{W}_{3L}) \ar[r]& 0
    }
    \]
    where rows are exact. As in the proof Proposition \ref{pi(W)}, $\theta_0$ and $\theta_1$ are identified with the map
    \[
      \Z\to\ell^\infty(\Z),\quad a\mapsto(\ldots,a,-a,a,-a,\ldots).
    \]
    Then it follows that the map $\pi_{2i}(\mathcal{W}_L)\to\pi_{2i}(\mathcal{W}_{3L})$ is an isomorphism for $1\le i\le L$.

    We shall show the map
    \begin{equation}
      \label{map W_L}
      \pi_{2i-1}(\mathcal{W}_L)\to\pi_{2i-1}(\mathcal{W}_{3L})
    \end{equation}
    is an isomorphism for $1\le i\le L$. By Lemma \ref{inclusion B}, we have
    \[
      \phi_1((a_j)_j)=(a_{3j}+a_{3j+1}+a_{3j+2})_j\quad\text{and}\quad\phi_2=\phi_1\times\phi_1.
    \]
    Then $\phi_1$ is surjective, and so is $\phi_2$. Thus \eqref{map W_L} is surjective. Now we suppose that $x\in\pi_{2i-1}(\mathcal{W}_L)$ is mapped to zero by \eqref{map W_L}. Then there is $(a,b)\in\ell^\infty(\Z)\times\ell^\infty(\Z)$ such that $(a,b)$ is mapped to $x$ and there is $c\in\ell^\infty(\Z)$ such that $\alpha_1(c)=\phi_2(a,b)$. Since $\phi_1$ is surjective, there is $\widetilde{c}\in\ell^\infty(\Z)$ such that $\phi_1(\tilde{c})=c$. Let $y=(a,b)-\alpha_0(\tilde{c})$. Then $y$ is mapped to $x$ and $\phi_2(y)=0$. Thus we may assume $\phi_2(a,b)=0$, or equivalently,
    \begin{equation}
      \label{(a,b)}
      a_{3j}+a_{3j+1}+a_{3j+2}=0\quad\text{and}\quad b_{3j}+b_{3j+1}+b_{3j+2}=0
    \end{equation}
    for each $j\in\Z$, where $a=(a_j)_j$ and $b=(b_j)_j$. The maps $\alpha_0$ and $\alpha_1$ are as in \eqref{alpha def}. Then by a straightforward computation, we can see that there is $d=(d_j)_j\in\ell^\infty(\Z)$ such that under the condition \eqref{(a,b)}, $\alpha_0(d)=(a,b)$, or equivalently,
    \[
      d_{2j}+d_{2j+1}=a_{j}\quad\text{and}\quad -d_{2j-1}-d_{2j}=b_{j}
    \]
    for each $j$. Indeed, if
    \begin{alignat*}{3}
      d_{6j}&=-b_{3j}&d_{6j+1}&=a_{3j}+b_{3j}\qquad&d_{6j+2}&=-a_{3j}+b_{3j+2}\\
      d_{6j+3}&=-a_{3j+2}-b_{3j+2}\qquad&d_{6j+4}&=a_{3j+2}&d_{6j+5}&=0,
    \end{alignat*}
    then $d\in\ell^\infty(\Z)$ and $\alpha_0(d)=(a,b)$. Thus $x=0$, implying the map \eqref{map W_L} is injective. Therefore, we obtain the map \eqref{map W_L} is an isomorphism for $1\le i\le L$.

    By Corollary \ref{pi(U^0)} together with the above computation, we obtain that for $i\ge 1$
    \[
      \pi_i(\mathcal{U}^0)=
      \begin{cases}
        \ell^\infty(\Z)_S & \text{$i$ is odd}\\
        \mathbb{Z} & \text{$i$ is even.}
      \end{cases}
    \]
    Since $\mathcal{U}$ is a topological group, all of its connected components have the same homotopy type. Then by Theorem \ref{index}, there is a homotopy equivalence
    \[
      \mathcal{U}\simeq\Z\times\mathcal{U}^0.
    \]
    Thus the proof is complete.
  \end{proof}


  \section{(End-)periodic finite propagation unitary operators}
  \label{section_periodic}

  In this section, we consider (end-)periodic finite propagation unitary operators. The space of periodic finite propagation unitary operators is of particular importance since it possesses the important homotopical information about $\mathcal{U}$ as we will see in Section 5. This section describes the homotopy types of the space of periodic finite propagation unitary operators $\widehat{\mathcal{U}}$ and the space of end-periodic finite propagation unitary operators $\widetilde{\mathcal{U}}$.

  The idea for describing the homotopy type of $\widehat{\mathcal{U}}$ is the same as the one in the previous section, that is, we will construct an approximating space for $\widehat{\mathcal{U}}$. The periodicity makes the approximating space drastically simple so that we can describe the homotopy type of $\widehat{\mathcal{U}}$, not only homotopy groups as in the previous section. We will give a product decomposition of $\widetilde{\mathcal{U}}$ involving $\widehat{\mathcal{U}}$, which describes its homotopy type.


  \subsection{Periodic finite propagation unitary operators}

  First, we define periodic finite propagation unitary operators. We say that a finite propagation unitary operator $U$ on $\ell^2(\C)$ is \emph{periodic} if
  \[
    S^nUS^{-n}=U
  \]
  for some $n>0$. Namely, $U=(U_{ij})$ is periodic if
  \[
  	U_{i+n,j+n}=U_{ij}
  \]
  for some $n>0$ and all $i,j$. For positive integers $L$ and $n$, we define
  \[
    \widehat{\mathcal{U}}_L(n)=\{ U\in\U(\B)\mid\prop(U)\le L,\,S^nUS^{-n}=U\}.
  \]
  We set $(L,n)\le(L',n')$ whenever $L\le L'$ and $n$ divides $n'$. Then we have the inclusion $\widehat{\mathcal{U}}_L(n)\to\widehat{\mathcal{U}}_{L'}(n')$ whenever $(L,n)\le(L',n')$, so the spaces $\widehat{\mathcal{U}}_L(n)$ form a direct system of spaces.

  \begin{definition}
    We define the space of periodic finite propagation unitary operators by
    \[
      \widehat{\mathcal{U}}=\lim_{\underset{(L,n)}{\longrightarrow}}\widehat{\mathcal{U}}_L(n)
    \]
    which is equipped with the weak topology with respect to $\widehat{\mathcal{U}}_L(n)$.
  \end{definition}

  Next, we define the spaces $\widehat{\mathcal{W}}_L$ which approximate $\widehat{\mathcal{U}}$. Let the group $\U(L)\times\U(L)$ act on the right on the space $\U(2L)\times\U(2L)$ by
  \[
    (U,U')\cdot(V,V')=(U(V\oplus V'),( (V')^{-1} \oplus V^{-1})U')
  \]
  for $U,U'\in\U(2L)$ and $V,V'\in\U(L)$, where for $A,B\in\U(L)$, $A\oplus B$ denotes the diagonal block matrix
  \[
    \begin{pmatrix}
      A&O\\O&B
    \end{pmatrix}\in\U(2L).
  \]
  We define the space $\widehat{\mathcal{W}}_L$ as the orbit space of this action. Since the action is free, we have the following.

  \begin{lemma}
    \label{fibration W period}
    There is a principal bundle
    \[
      \U(L)\times\U(L)\to\U(2L)\times\U(2L)\to\widehat{\mathcal{W}}_L.
    \]
  \end{lemma}

  Let $\Delta_k\colon\U(L)\to\U(B_k(L))$ denote the inclusion such that
  \[
    \Delta_k(U)=
    \begin{pmatrix}
      \ddots \\
      & U \\
      &&  U \\
      &&  & U & \\
      &&&& \ddots
    \end{pmatrix}
  \]
  for $U\in\U(L)$, where the $(1,1)$-entry of each $U$ in $\Delta_k(U)$ is located at the $(k+nL,k+nL)$-entry of the whole matrix for some integer $n$. Note that $\Delta_k(U)$ is periodic such that
  \[
    S^L\Delta_k(U) S^{-L} = \Delta_k(U).
  \]

  By the definition of $\widehat{\mathcal{W}}_L$, the map
  \[
    \phi_L\colon\widehat{\mathcal{W}}_L\to\mathcal{U}^0,\quad[U,U']\mapsto\Delta_0(U)\Delta_{-L}(U')
  \]
  is a well-defined homeomorphism onto its image. For $L\mid L'$, there is a map $\alpha_{L',L}\colon\U(2L)\times\U(2L)\to\U(2L')\times\U(2L')$ given as
  \[
  \alpha_{L',L}(U,U')=\left(
  \begin{pmatrix}
    U \\
  & \ddots \\
  & & U
  \end{pmatrix}
  \begin{pmatrix}
  I_L \\
  & U' \\
  & & \ddots \\
  & & & U' \\
  & & & & I_L
  \end{pmatrix}
  ,
  \begin{pmatrix}
  I_{L'-L} \\
  & U' \\
  & & I_{L'-L}
  \end{pmatrix}
  \right)
  \]
  where $I_L$ denotes the identity matrix. We can see that $\alpha_{L',L}$ induces a map
  \[
    \widehat{\alpha}_{L',L}\colon\widehat{\mathcal{W}}_L\to\widehat{\mathcal{W}}_{L'},\quad[U,U']\mapsto[\alpha_{L',L}(U,U')]
  \]
  such that
  \begin{equation}
    \label{phi_L}
    \phi_{L'}\circ\widehat{\alpha}_{L',L}=\phi_L.
  \end{equation}
  For $L\mid L'$ and $L'\mid L''$,
  \[
    \phi_{L''}\circ\widehat{\alpha}_{L'',L'}\circ\widehat{\alpha}_{L',L}
    =
    \phi_{L}
    =
    \phi_{L''}\circ\widehat{\alpha}_{L'',L}.
  \]
  Then since $\phi_{L''}$ is injective, we get $\widehat{\alpha}_{L'',L'}\circ\widehat{\alpha}_{L',L}=\widehat{\alpha}_{L'',L}$, and so the spaces $\widehat{\mathcal{W}}_L$ form a direct system of spaces, and we define
  \[
    \widehat{\mathcal{W}}=\lim_{\underset{L}{\longrightarrow}}\widehat{\mathcal{W}}_L.
  \]
  Let $\widehat{\mathcal{U}}^0$ denote the identity component of $\widehat{\mathcal{U}}$.

  \begin{lemma}
    \label{U^0 period}
  There is a homeomorphism
  \[
  	\widehat{\mathcal{U}}^0\cong\widehat{\mathcal{W}}.
  \]
\end{lemma}

\begin{proof}
  By \eqref{phi_L}, we get a map
  \[
    \phi\colon\widehat{\mathcal{W}}\to\mathcal{U}^0.
  \]
  Since the map $\widehat{\mathcal{W}}_L\to\mathcal{U}^0$ is a homeomorphism into, the map $\phi$ is also a homeomorphism into. Then it suffices to show $\phi(\widehat{\mathcal{W}})=\widehat{\mathcal{U}}^0$. Clearly, every element of $\phi(\widehat{\mathcal{W}})$ is periodic, so $\phi(\widehat{\mathcal{W}})\subset\widehat{\mathcal{U}}^0$. As in Remark \ref{rem_decomp_periodic}, for each $U\in\widehat{\mathcal{U}}_L(n)\cap\mathcal{U}^0$ and an integer $M$ divisible by $L$ and $n$, there are $V,W\in\U(2M)$ such that
  \[
    U=\Delta_0(V)\Delta_{-M}(W).
  \]
  Then $\phi(\widehat{\mathcal{W}})\supset\widehat{\mathcal{U}}^0$, so $\phi(\widehat{\mathcal{W}})=\widehat{\mathcal{U}}^0$. Thus the proof is done.
\end{proof}


  \subsection{Cohomology of $\widehat{\mathcal{W}}_L$}

  We compute the cohomology of $\widehat{\mathcal{W}}_L$, which will be used to give the homotopy decomposition of $\widehat{\mathcal{U}}$. Recall that the cohomology of $B\U(n)$ and $\U(n)$ are given by
  \begin{align*}
    H^*(B \U(n))&=\mathbb{Z}[c_1,c_2,\ldots,c_n]\\
    H^*(\U(n))&=\Lambda(e_1,e_3,\ldots,e_{2n-1})
  \end{align*}
  where $c_i$ is the $i$-th Chern class of the universal bundle and $e_i$ is characterized by $\tau(e_{2i-1})=c_i$ for the transgression $\tau$ in the Serre spectral sequence of the universal bundle $\U(n)\to E\U(n)\to B\U(n)$. Let $\Gr_n(\mathbb{C}^{2n})$ denote the Grassmannian of $n$-dimensional subspaces in $\mathbb{C}^{2n}$. Let $G\to E\to B$ be a principal bundle, where $B$ is a CW complex. Then there is a pullback diagram
  \[
    \xymatrix{
      E\ar[r]\ar[d]&EG\ar[d]\\
      B\ar[r]&BG.
    }
  \]
  Since $EG$ is contractible, $E$ is homotopy equivalent to the homotopy fiber of the map $B\to BG$. Then we get a homotopy fibration $E\to B\to BG$. We will use freely this fact. Since there is a principal bundle $\U(n)\times\U(n)\to\U(2n)\to\Gr_n(\mathbb{C}^{2n})$, there is a homotopy fibration
  \[
    \U(2n)\to\Gr_n(\mathbb{C}^{2n})\to B\U(n)\times B\U(n).
  \]

  \begin{lemma}
    \label{Cartan}
  	In the Serre spectral sequence of the homotopy fibration above, we have
    \[
      \tau(e_{2k-1})=\sum_{i+j=k}c_i\times c_j.
    \]
  \end{lemma}

  \begin{proof}
  	There is a homotopy commutative diagram
  	\[
      \xymatrix{
        \U(2n)\ar[r]\ar@{=}[d]&\Gr_n(\mathbb{C}^{2n})\ar[r]\ar[d]&B\U(n)\times B\U(n)\ar[d]^g\\
        \U(2n)\ar[r]&E\U(2n)\ar[r]&B\U(2n)
      }
    \]
  	where rows are homotopy fibrations and the map $g$ is induced from the diagonal block inclusion $\U(n)\times\U(n)\to\U(2n)$. Compare the Serre spectral sequences of the top and the bottom homotopy fibrations. Then since $g^*(c_k)=\sum_{i+j=k}c_i\times c_j$, we get
    \[
      \tau(e_{2k-1})=g^*(c_k)=\sum_{i+j=k}c_i\times c_j
    \]
    as stated.
  \end{proof}

  From Lemma \ref{fibration W period}, the quotient of the first projection $\U(2L)\times\U(2L)\to\U(2L)$ by the action of $\U(L)\times\U(L)$ is a map
  \[
  \rho\colon\widehat{\mathcal{W}}_L\to\Gr_{L}(\mathbb{C}^{2L}),
  \]
  which is the projection of a fiber bundle with fiber $\U(2L)$.

  \begin{proposition}
    \label{W cohomology}
  	The rational cohomology Serre spectral sequence of the fiber bundle
  	\[
  	\U(2L)\to\widehat{\mathcal{W}}_L\xrightarrow{\rho}\Gr_{L}(\mathbb{C}^{2L})
  	\]
  	collapses.
  \end{proposition}

  \begin{proof}
    By the definition of $\widehat{\mathcal{W}}_L$, there is a homotopy commutative diagram
    \[
  	  \xymatrix{
  	    \U(2L)\times\U(2L)\ar[r]\ar[d]_{1\times\iota}&\widehat{\mathcal{W}}_L\ar[r]\ar[d]&B\U(L)\times B\U(L)\ar[d]\\
  	    \U(2L)\times\U(2L)\ar[r]&E\U(2L)\times E\U(2L)\ar[r]&B\U(2L)\times B\U(2L)
      }
    \]
    where rows are homotopy fibrations and $\iota\colon\U(2L)\to\U(2L)$ is given by $\iota(A)=A^{-1}$ for $A\in\U(2L)$. The maps between the rows are obtained from the fact that the map $\U(L)\times\U(L)\to\U(2L)\times\U(2L)$ in Lemma \ref{fibration W period} is a homomorphism into the first component and an anti-homomorphism into the second component.  Consider the Serre spectral sequence of the top row. The $E_2$-term is given as
    \[
      E_2^{\ast,\ast}\cong
      H^\ast(B\U(L)\times B\U(L);\mathbb{Q})
      	\otimes H^\ast(\U(2L)\times\U(2L);\mathbb{Q}).
    \]
    Comparing with the bottom row, the first non-trivial differentials on the generators of $E_2^{0,\ast}$ are the transgressions
    \[
      \tau(e_{2k-1}\times 1)=\sum_{i+j=k}c_i\times c_j,\quad\tau(1\times e_{2k-1})=-\sum_{i+j=k}c_i\times c_j.
    \]
    This implies that the element $e_{2k-1}\times1+1\times e_{2k-1}$ ($1\le k\le 2L$) is contained in $E_\infty^{0,\ast}$.
    Then the composite
    \begin{align}
    \U(2L)\xrightarrow{\text{the second inclusion}}
    \U(2L)\times\U(2L)\to\widehat{\mathcal{W}}_L
    \tag{$\ast$}
    \end{align}
    induces a surjection on cohomology.
	Since the fiber inclusion of the fiber bundle
    \[
    \U(2L)\to\widehat{\mathcal{W}}_L\xrightarrow{\rho}\Gr_{L}(\mathbb{C}^{2L})
    \]
    coincides with the composite ($\ast$), we can apply the Leray--Hirsch theorem.
    Thus the spectral sequence collapses.
  \end{proof}


  \subsection{Homotopy decomposition of $\widehat{\mathcal{U}}$}

  We give a homotopy decomposition of $\widehat{\mathcal{U}}$. We will use the following lemma to prove it.

  \begin{lemma}
  	\label{colim fibration}
  	Let $\{F_\lambda\to E_\lambda\to B_\lambda\}_{\lambda\in\Lambda}$ be a direct system of homotopy fibrations such that all $F_\lambda,E_\lambda,B_\lambda$ are path-connected and all maps $F_\lambda\to F_\mu,E_\lambda\to E_\mu,B_\lambda\to B_\mu$ are cofibrations. Then
  	$$\lim_{\longrightarrow}F_\lambda\to\lim_{\longrightarrow}E_\lambda\to\lim_{\longrightarrow}B_\lambda$$
  	is a homotopy fibration.
  \end{lemma}

  \begin{proof}
  	Since direct limits commute with exact sequences, there is an exact sequence
  	$$\cdots\to\lim_{\longrightarrow}\pi_*(F_\lambda)\to\lim_{\longrightarrow}\pi_*(E_\lambda)\to\lim_{\longrightarrow}\pi_*(B_\lambda)\to\cdots.$$
  	Since all maps $F_\lambda\to F_\mu,E_\lambda\to E_\mu,B_\lambda\to B_\mu$ are cofibrations,
  	$$\lim_{\longrightarrow}\pi_*(F_\lambda)\cong\pi_*(\lim_{\longrightarrow}F_\lambda),\quad\lim_{\longrightarrow}\pi_*(E_\lambda)\cong\pi_*(\lim_{\longrightarrow}E_\lambda),\quad\lim_{\longrightarrow}\pi_*(B_\lambda)\cong\pi_*(\lim_{\longrightarrow}B_\lambda).$$
  	Thus the proof is done.
  \end{proof}
  \begin{remark}
  This lemma has several variants.
  For example, Lemma 7.4.1 of \cite{Hov} states that a sequential colimit of fibrations is again a fibration while we prove it up to homotopy equivalence in some sense.
  \end{remark}

  \begin{proposition}
    \label{fibration W period 2}
    There is a principal homotopy fibration
    \[
  		\overset{\circ}{\prod_{n\ge 1}} K(\Q,2n-1)\to\widehat{\mathcal{W}}\to B\U(\infty).
  	\]
  \end{proposition}

  \begin{proof}
    Note that the inclusion $\widehat{\mathcal{W}}_L\to\widehat{\mathcal{W}}_{mL}$ is $\U(2L)$-equivariant along the inclusion
  	\[
  		\U(2L)\to\U(2mL),
  		\quad
  		V\to V^{\oplus m}.
  	\]
  	It induces the following commutative diagram.
  	\[
  		\xymatrix{
  			\U(2L)\ar[r]\ar[d]&\widehat{\mathcal{W}}_L\ar[r]^-{\rho}\ar[d]&\Gr_{L}(\mathbb{C}^{2L})\ar[d]\\
  			\U(2mL)\ar[r]&\widehat{\mathcal{W}}_{mL}\ar[r]^-{\rho}&\Gr_{mL}(\mathbb{C}^{2mL})
  		}
  	\]
  	Clearly, the direct limit of the sequence of the natural inclusions
  	$$\cdots\to\Gr_n(\mathbb{C}^{2n})\to\Gr_{n+1}(\mathbb{C}^{2n+2})\to\cdots$$
  	is $B\U(\infty)$. Note that the direct limit of the sequence
  	$$\Z\xrightarrow{n}\Z\xrightarrow{n+1}\Z\xrightarrow{n+2}\cdots$$
  	is $\Q$. The diagonal block inclusion $\U(n)\to\U(mn)$ is homotopic to the composition of the $m$-power map $\U(n)\to\U(n)$ and the canonical inclusion $\U(n)\to\U(mn)$. Therefore the diagonal block inclusion $\U(n)\to\U(mn)$ induces multiplication by $m$ on the homotopy groups in the stable range. Thus the direct limit of the sequence of diagonal block inclusions
    \[
      \U(n)\to\U(kn)\to\U(k(k+1)n)\to\cdots.
    \]
    is the rationalization of $\U(\infty)$. It is well known that the rationalization of $\U(\infty)$ is homotopy equivalent to $\overset{\circ}{\prod}_{n\ge 1} K(\Q,2n-1)$, and so by Lemma \ref{colim fibration}, we obtain the homotopy fibration in the statement, which is obviously principal.
  \end{proof}

  Now we prove the homotopy decomposition of $\widehat{\mathcal{U}}$.

  \begin{theorem}
  	\label{decomp period}
  	There is a homotopy equivalence
    \[
      \widehat{\mathcal{U}}\simeq\Z\times B\U(\infty)\times\overset{\circ}{\prod_{n\ge 1}} K(\Q,2n-1).
    \]
  \end{theorem}

  \begin{proof}
    We show that the principal fibration of Proposition \ref{fibration W period 2} is trivial.
    Since the natural map $H^*(\mathrm{Gr}_{L+1}(\C^{2L+2});\mathbb{Q})\to H^*(\mathrm{Gr}_L(\C^{2L});\mathbb{Q})$ is surjective,
	\[
	  \underset{L\to\infty}{\mathrm{lim}^1}H^*(\mathrm{Gr}_L(\C^{2L});\mathbb{Q})=0.
	\]
	Then by Proposition \ref{W cohomology} and the proof of Proposition \ref{fibration W period 2}, its rational cohomology Serre spectral sequence collapses.
	Consider the map of spectral sequences induced from the map to the universal fibration
	\[
	\xymatrix{
		\displaystyle{\overset{\circ}{\prod_{n\ge 1}}K(\Q,2n-1)} \ar[r] \ar@{=}[d]
		& \widehat{\mathcal{W}} \ar[r] \ar[d]
		& B\U(\infty) \ar[d] \\
		\displaystyle{\overset{\circ}{\prod_{n\ge 1}}K(\Q,2n-1)} \ar[r]
		& \ast \ar[r]
		& \displaystyle{\overset{\circ}{\prod_{n\ge 1}}K(\Q,2n)}
	}
	\]
	Since the fundamental classes in $H^{2n}(K(\mathbb{Q},2n);\mathbb{Q})$ for $n\ge1$ are the transgressions from $H^{2n-1}(K(\mathbb{Q},2n-1);\mathbb{Q})$ and the spectral sequence of the top row collapses,
	the classifying map
    \[
      B\U(\infty)\to\overset{\circ}{\prod_{n\ge 1}} K(\Q,2n).
    \]
	must be null-homotopic.
	Thus the principal fibration of Proposition \ref{fibration W period 2} is trivial and we have the homotopy equivalence
	\[
	\widehat{\mathcal{W}}
	\simeq
	B\U(\infty)\times\displaystyle{\overset{\circ}{\prod_{n\ge 1}}K(\Q,2n-1)}
	\]
    If $U\in\widehat{\mathcal{U}}^0$, then $S^nU\in\widehat{\mathcal{U}}$ for any integer $n$. Since $\ind(S^nU)=n\ind(S)+\ind(U)=n$ by Theorem \ref{index}, we get $\pi_0(\widehat{\mathcal{U}})\cong\Z$. Then as in the proof of Theorem \ref{main 1}, we obtain
    \[
      \widehat{\mathcal{U}}\simeq\Z\times\widehat{\mathcal{U}}^0.
    \]
    Thus, by the homotopy equivalence $\widehat{\mathcal{U}}^0\simeq\widehat{\mathcal{W}}$ in Lemma \ref{U^0 period}, the proof is complete.
  \end{proof}

  \begin{corollary}
    \label{pi(U) period}
    There is an isomorphism
    \[
      \pi_i(\widehat{\mathcal{U}})\cong
      \begin{cases}
        \mathbb{Q}&i\text{ is odd}\\
        \mathbb{Z}&i\text{ is even.}
      \end{cases}
    \]
  \end{corollary}


  \subsection{End-periodic finite propagation unitary operators}

  First, we define end-periodic finite propagation unitary operators. We say that a finite propagation unitary operator $U$ on $\ell^2(\C)$ is \emph{end-periodic} if $S^nUS^{-n}-U$ has only finitely many non-zero entries for some $n$. For positive integers $L, m,n$, let $\widetilde{\mathcal{U}}_L(m,n)$ denote the set of end-periodic finite propagation unitary operators $U$ on $\ell^2(\C)$ such that $\prop(U)\le L$ and $S^nUS^{-n}-U$ is non-trivial only in the central $2m\times2m$ region of the matrix. Clearly, periodic finite propagation unitary operators are end-periodic and so
  \[
    \widehat{\mathcal{U}}_L(n)\subset\widetilde{\mathcal{U}}_L(m,n).
  \]

  If $L\le L'$, $n$ divides $n'$ and $m\le m'$, there is the inclusion
  \[
    \widetilde{\mathcal{U}}_L(m,n)\to\widetilde{\mathcal{U}}_{L'}(m',n')
  \]
  and so the spaces $\widetilde{\mathcal{U}}_L(m,n)$ form a direct system of spaces.

  \begin{definition}
    The space of end-periodic finite propagation unitary operators is defined by
    \[
      \widetilde{\mathcal{U}}=\lim_{\underset{(L,m,n)}{\longrightarrow}}\widetilde{\mathcal{U}}_L(m,n)
    \]
    which is given the weak topology with respect to $\widetilde{\mathcal{U}}_L(m,n)$.
  \end{definition}

  We can embed $U(2m)$ into $\mathcal{U}$ by sending a unitary matrix to the finite propagation unitary operator
  \[
    \begin{pmatrix}
      \ddots\\
      &1\\
      &&U\\
      &&&1\\
      &&&&\ddots
    \end{pmatrix}
  \]
  where the $(1,1)$-entry of $U$ is at the $(-m,m)$-entry of the whole matrix. Then we may assume that $\U(2m)$ is a subgroup of $\mathcal{U}$.

  \begin{lemma}
  \label{decomp end-period L}
    The map
    \[
      \U(m)\times\widehat{\mathcal{U}}_L(n)\to\widetilde{\mathcal{U}}_L(m,n),\quad(U,V)\mapsto UV
    \]
    is a homeomorphism.
  \end{lemma}

  \begin{proof}
    For any $U\in\widetilde{\mathcal{U}}_L(m,n)$, there is a unique periodic finite propagation unitary operator $r(U)\in\widehat{\mathcal{U}}_L(n)$ which coincides with $U$ outside the central $2m\times 2m$ region of the matrix $U$. Then we get a map
    \[
      r\colon\widetilde{\mathcal{U}}_L(m,n)\to\widehat{\mathcal{U}}_L(n)
    \]
    which is a retraction of the inclusion $\widehat{\mathcal{U}}_L(n)\to\widetilde{\mathcal{U}}_L(m,n)$. By the definition of $r$, the map
    \[
      \widetilde{\mathcal{U}}_L(m,n)\to\U(m)\times\widehat{\mathcal{U}}_L(n),\quad U\mapsto(r(U)^{-1}U,r(U))
    \]
    is the inverse of the map in the statement, completing the proof.
  \end{proof}

  Now we are ready to give the product decomposition of $\widetilde{\mathcal{U}}$.

  \begin{theorem}
    \label{decomp end-period}
    There is a homeomorphism
    \[
      \widetilde{\mathcal{U}}\cong\U(\infty)\times\widehat{\mathcal{U}}.
    \]
  \end{theorem}

  \begin{proof}
    By definition, we have
    \[
      \widetilde{\mathcal{U}}=\lim_{m\to\infty}\lim_{\underset{(L,n)}{\longrightarrow}}\widetilde{\mathcal{U}}_L(m,n).
    \]
    Since $U(m)$ is compact, taking a product with $U(m)$ commutes with taking a direct limit. Then since the homeomorphism of Lemma \ref{decomp end-period L} is natural with respect to $(L,n)$, we get
    \[
      \lim_{\underset{(L,n)}{\longrightarrow}}\widetilde{\mathcal{U}}_L(m,n)\cong\lim_{\underset{(L,n)}{\longrightarrow}}(\U(m)\times\widehat{\mathcal{U}}_L(n))\cong\U(m)\times\lim_{\underset{(L,n)}{\longrightarrow}}\widehat{\mathcal{U}}_L(n)=\U(m)\times\widehat{\mathcal{U}}.
    \]
    Since $\widehat{\mathcal{U}}$ is a compactly generated Hausdorff space, taking a product with $\widehat{\mathcal{U}}$ also commutes with taking a direct limit. Then
    \[
      \lim_{m\to\infty}(\U(m)\times\widehat{\mathcal{U}})\cong(\lim_{m\to\infty}\U(m))\times\widehat{\mathcal{U}}=\U(\infty)\times\widehat{\mathcal{U}}.
    \]
    Thus the proof is complete.
  \end{proof}

  By Theorems \ref{decomp period} and \ref{decomp end-period}, we obtain the following corollary.

  \begin{corollary}
    There is a homotopy equivalence
    \[
      \widetilde{\mathcal{U}}\simeq\Z\times\U(\infty)\times B\U(\infty)\times\overset{\circ}{\prod_{n\ge 1}} K(\Q,2n-1).
    \]
  \end{corollary}


  \section{Proof of Theorem \ref{main 2}}

  In this section, we prove Theorem \ref{main 2}. There are three steps in the proof. The first step is to show that the abelian group $\ell^\infty(\Z)_S$ is a $\mathbb{Q}$-vector space, and the second step is to describe the homotopy type of a rational topological monoid having trivial even homotopy groups. The third step is to show that the principal bundle in Lemma \ref{fibration W} is trivial by using the results in the former steps together with Theorem \ref{decomp period}.


  \subsection{The abelian group $\ell^\infty(\Z)_S$}

  We aim to prove that the abelian group $\ell^\infty(\Z)_S$ is a $\mathbb{Q}$-vector space. This is unexpected because it is obvious that $\ell^\infty(\Z)$ is not a $\mathbb{Q}$-vector space. Recall that an abelian group $A$ is called divisible if for any $a\in A$ and a positive integer $n$, there is $b\in B$ such that $a=nb$. If $b$ is unique for every $a$ and $n$, then $A$ is called uniquely divisible. It is clear that an abelian group is a $\mathbb{Q}$-vector space if and only if it is uniquely divisible.

  \begin{lemma}
    \label{divisible}
    The abelian group $\ell^\infty(\mathbb{Z})_S$ is divisible.
  \end{lemma}

  \begin{proof}
  Take any $a=(a_j)_j\in\ell^\infty(\mathbb{Z})$ and a positive integer $n$.
  Then we need to find $b=(b_j)_j,c=(c_j)_j\in\ell^\infty(\mathbb{Z})$ satisfying
  \[
  a-nb=(1-S)(c).
  \]
  We construct such $b$ and $c$ by induction. Set $c_0=0$. Suppose that $c_j$ and $c_{-j}$ are given for some $j\ge0$ such that $0\le c_j<n$ and $0\le c_{-j}<n$. Then there is $b_j$ satisfying $0\le-a_j+nb_j+c_j<n$. We set $c_{j+1}=-a_j+nb_j+c_j$. There is also $b_{-j-1}$ satisfying $0\le a_{-j-1}-nb_{-j-1}+c_{-j}<n$. Then we also set $c_{-j-1}=a_{-j-1}-nb_{-j-1}+c_{-j}$. By definition, $0\le c_j<n$ for all $j$, so $c\in\ell^\infty(\Z)$. On the other hand, we have
  \[
    |b_j|=\frac{1}{n}|a_j-c_j+c_{j+1}|<\frac{|a_j|}{n}+2\quad\text{and}\quad|b_{-j-1}|=\frac{1}{n}|a_{-j-1}-c_{-j-1}+c_j|<\frac{|a_{-j-1}|}{n}+2,
  \]
  implying $b\in\ell^\infty(\Z)$.
  Also, $a-nb=(1-S)(c)$ holds by construction, completing the proof.
  \end{proof}

  \begin{lemma}
    \label{torsion-free}
    The abelian group $\ell^\infty(\mathbb{Z})_S$ is torsion-free.
  \end{lemma}

  \begin{proof}
    Suppose that the equality
    \[
      na=(1-S)(c)
    \]
    holds for some $a=(a_j)_j,c=(c_j)_j\in\ell^\infty(\mathbb{Z})$ and a positive integer $n$. Since the constant sequence $b$ satisfies $(1-S)(b)=0$, we may assume $c_0=0$. Then by induction, we can see that $c_j$ is divisible by $n$ for all $j$, implying there is $c'\in\ell^\infty(\mathbb{Z})$ satisfying $c=nc'$. Thus
    \[
      n(a-(1-S)(c'))=0.
    \]
    Clearly, $\ell^\infty(\mathbb{Z})$ is torsion-free, so we obtain
    \[
      a=(1-S)c'.
    \]
    Thus we have proved that if $x\in\ell^\infty(\Z)_S$ satisfies $nx=0$ for some positive integer $n$, then $x=0$. Namely, $\ell^\infty(\Z)_S$ is torsion free.
  \end{proof}

  Now we are ready to prove the following proposition.

  \begin{proposition}
    \label{coinvariant}
    The abelian group $\ell^\infty(\mathbb{Z})_S$ is a $\mathbb{Q}$-vector space.
  \end{proposition}

  \begin{proof}
    As mentioned above, we need to prove $\ell^\infty(\Z)_S$ is uniquely divisible. By Lemma \ref{divisible}, $\ell^\infty(\Z)_S$ is divisible, that is, for any $a\in\ell^\infty(\Z)_S$ and any positive integer $n$, there is $b\in\ell^\infty(\Z)_S$ satisfying
    \[
      a=nb.
    \]
    It remains to show that such an element $b$ is unique. Suppose that $b'\in\ell^\infty(\Z)_S$ satisfies $a=nb'$. Then $n(b-b')=0$. It follows from Lemma \ref{torsion-free} that $b=b'$, completing the proof.
  \end{proof}


  \subsection{Homotopy type of a rational H-space}

  It is well known that a rational H-space of finite rational type is of the homotopy type of the product of Eilenberg-MacLane spaces \cite[Example 3, p. 143]{FHT}. Here we consider the homotopy type of a rational H-space which is not of finite rational type.

  \begin{lemma}
    \label{rational H-space}
    Let $X$ be a connected homotopy associative H-space which is a CW-complex. Suppose that the even homotopy groups of $X$ are trivial and the odd homotopy groups of $X$ are $\mathbb{Q}$-vector spaces. Then there is a homotopy equivalence
    \[
      X\simeq\prod_{i=1}^{\circ}K(\pi_{2i-1}(X),2i-1).
    \]
  \end{lemma}

  \begin{proof}
    Let $B_i$ be a basis of $\pi_{2i-1}(X)$. Then there is a direct system of finite sets $\{B_i^\lambda\}_{\lambda\in\Lambda_i}$ such that
    \[
      B_i=\lim_{\longrightarrow}B_i^\lambda.
    \]
    We choose a representative of $b\in B_i\subset\pi_{2i-1}(X)$, and denote it by the same symbol $b$. For $b\in B_i$, let $S(b)$ denote the rational $(2i-1)$-sphere. We give a total ordering on $\Lambda_i$ and consider the map
    \[
      \phi_i^\lambda\colon\prod_{b\in B_i^\lambda}S(b)\to X,\quad(x_b)_{b\in B_i^\lambda}\mapsto \prod_{b\in B_i^\lambda}b(x_b)
    \]
    where the product is taken following the order of $\lambda_i$. Then $\phi_i^\lambda$ is natural with respect to $\lambda$, so we get a map
    \[
      \phi_i\colon\lim_{\underset{\lambda}{\longrightarrow}}\prod_{b\in B_i^\lambda}S(b)\to X.
    \]
    Since $S(b)\simeq K(\mathbb{Q},2i-1)$ and $\pi_{2i-1}(X)=\bigoplus_{b\in B_i}\mathbb{Q}\langle b\rangle$, we have
    \[
      \lim_{\underset{\lambda}{\longrightarrow}}\prod_{b\in B_i^\lambda}S(b)\simeq K(\pi_{2i-1}(X),2i-1).
    \]
    Then we get a map $\bar{\phi}_n\colon\prod_{i=1}^nK(\pi_{2i-1}(X),2i-1)\to X$ such that
    \[
      \bar{\phi}_n(x_1,\ldots,x_n)=(\cdots((\phi_1(x_1)\phi_2(x_2))\phi_3(x_3))\cdots\phi_{n-1}(x_{n-1}))\phi_n(x_n)
    \]
    which is natural with respect to $n$. Thus we can take the direct limit with respect to $n$, and obtain a map
    \[
      \bar{\phi}\colon\prod_{i=1}^{\circ}K(\pi_{2i-1}(X),2i-1)\to X.
    \]
    By the construction, $\bar{\phi}$ is an isomorphism in homotopy groups. Now both $\prod_{i=1}^{\circ}K(\pi_{2i-1}(X),2i-1)$ and $X$ are CW-complexes, the proof is complete by the J.H.C. Whitehead theorem.
  \end{proof}


  \subsection{Proof of Theorem \ref{main 2}}

  Consider the projection $\pi\colon\U(B_{-L}(2L))\to\U(2L)$ onto the block unitary matrix containing the $(0,0)$-entry. Then there is a commutative diagram
  \[
    \xymatrix{
      \U(B_0(L))\ar[r]\ar[d]&\U(B_0(2L))\times\U(B_{-L}(2L))\ar[r]\ar[d]&\mathcal{W}_L\ar[d]\\
      \U(L)\times\U(L)\ar[r]&\U(2L)\ar[r]&\Gr_L(\C^{2L})
    }
  \]
  where the top row is the principal fiber bundle of Lemma \ref{fibration W}, the bottom row is the canonical principal bundle and the central map $\U(B_0(2L))\times\U(B_{-L}(2L))\to\U(2L)$ is given by $(U,(V_j)_j)\mapsto V_0$, where $V_0$ means the $0$-th block as in Section \ref{section_index}.

  \begin{lemma}
    \label{pi(W) even}
    The map $\mathcal{W}_L\to\Gr_L(\C^{2L})$ is an isomorphism in $\pi_{2i}$ for $i\le L$.
  \end{lemma}

  \begin{proof}
    By the commutative diagram above, we get a commutative diagram
    \[
      \xymatrix{
        \pi_{2i}(\mathcal{W}_L)\ar[r]\ar[d]&\pi_{2i-1}(\U(B_0(L)))\ar[d]\\
        \pi_{2i}(\Gr_L(\C^{2L}))\ar[r]&\pi_{2i-1}(\U(L))\times\pi_{2i-1}(\U(L)).
      }
    \]
    For $i\le L$, the top map is identified with the map
    \[
      \Z\to\ell^\infty(\Z),\quad a\mapsto(\ldots,a,-a,a,-a,\ldots)
    \]
    as in the proof of Theorem \ref{main 1}, and the bottom map is identified with the map
    \[
      \Z\to\Z\times\Z,\quad a\mapsto(a,-a).
    \]
    Moreover, since the map $\U(B_0(L))\to\U(L)\times\U(L)$ is the restriction of the projection $\U(B_{-L}(2L))\to\U(2L)$, the right vertical map is identified with the map
    \[
      \ell^\infty(\Z)\to\Z\times\Z,\quad(a_j)_j\mapsto(a_0,a_1)
    \]
    for $i\le L$. Thus we can see that the left vertical map is an isomorphism.
  \end{proof}

  By the construction, the map $\mathcal{W}_L\to\Gr_L(\mathbb{C}^{2L})$ is natural with respect to $L$, so we get a map
  \[
    \mathcal{W}\to B\U(\infty)
  \]
  where $\lim_{L\to\infty}\Gr_L(\mathbb{C}^{2L})= B\U(\infty)$. By Lemma \ref{pi(W) even} and the proof of Theorem \ref{main 1}, we have the following lemma.

  \begin{lemma}
    \label{pi(U) even}
    The map $\mathcal{W}\to B\U(\infty)$ is an isomorphism in the even homotopy groups.
  \end{lemma}

  Consider the inclusion
  \[
    \mu_L\colon\U(L)\times\U(L)\to\U(2L),\quad(A,B)\mapsto\begin{pmatrix}A&O\\O&B\end{pmatrix}.
  \]
  Then this induces a map
  \[
    \bar{\mu}_L\colon \U(B_k(L))\times\U(B_k(L))\to\U(B_k(2L))
  \]
  which is compatible with the projection $\pi\colon\U(B_{-L}(2L))\to\U(2L)$. Then we get a commutative diagram
  \[
    \xymatrix{
      \mathcal{W}_L\times\mathcal{W}_L\ar[r]\ar[d]&\mathcal{W}_{2L}\ar[d]\\
      \Gr_L(\mathbb{C}^{2L})\times\Gr_L(\mathbb{C}^{2L})\ar[r]&\Gr_{2L}(\mathbb{C}^{4L})
    }
  \]
  which is natural with respect to $L$, up to homotopy, where the bottom map is the canonical inclusion. Clearly, $\mathcal{W}$ is an H-space by the induced map $\mathcal{W}\times\mathcal{W}\to\mathcal{W}$. (Actually, $\mathcal{W}$ is an $A_\infty$-space by this product.) Moreover, since the canonical inclusion $\Gr_L(\mathbb{C}^{2L})\times\Gr_L(\mathbb{C}^{2L})\to\Gr_{2L}(\mathbb{C}^{4L})$ induces the standard multiplication of $B\U(\infty)$, the map $\mathcal{W}\to B\U(\infty)$ is an H-map. (This map is actually an $A_\infty$-map.) Let $\mathcal{F}$ denote the homotopy fiber of this H-map.

  \begin{lemma}
    \label{deloop}
    The homotopy fibration $\mathcal{F}\to\mathcal{W}\to B\U(\infty)$ is an H-fibration, that is, all spaces are H-spaces and all maps are H-maps.
  \end{lemma}

  \begin{proof}
    Since the homotopy fiber of an H-map is an H-space and the canonical map from the homotopy fiber is an H-map, the statement follows from the above observation that the map $\mathcal{W}\to B\U(\infty)$ is an H-map.
  \end{proof}

  Now we determine the homotopy type of $\mathcal{F}$.

  \begin{lemma}
    \label{F}
    There is a homotopy equivalence
    \[
      \mathcal{F}\simeq\overset{\circ}{\prod_{i\ge 1}} K(\ell^\infty(\Z)_S,2i-1).
    \]
  \end{lemma}

  \begin{proof}
    By Corollary \ref{pi(U^0)} and Lemma \ref{pi(U) even}, there is an isomorphism for $i\ge 1$
    \[
      \pi_i(\mathcal{F})\cong
      \begin{cases}
        \ell^\infty(\Z)_S&i\text{ is odd}\\
        0&i\text{ is even.}
      \end{cases}
    \]
    Then since loop spaces are equivalent to topological monoids, the homotopy equivalence in the statement is obtained by Lemma \ref{rational H-space}.
  \end{proof}

  We will use the following simple lemma, which is well known in homotopy theory.

  \begin{lemma}
    \label{fibration trivial}
    An H-fibration $F\xrightarrow{i}E\xrightarrow{p}B$ is trivial if and only if the map $p$ has a homotopy section, where $F,E,B$ are CW-complexes.
  \end{lemma}

  \begin{proof}
    The only if part is obvious, so we only prove the if part. Let $s\colon B\to E$ be a homotopy section of $p$. Define a map
    \[
      \phi\colon B\times F\to E,\quad (x,y)\mapsto s(x)i(y).
    \]
    Then since $p$ is an H-map, there is a homotopy commutative diagram
    \[
      \xymatrix{
        F\ar[r]\ar@{=}[d]&B\times F\ar[r]\ar[d]^\phi&B\ar@{=}[d]\\
        F\ar[r]^i&E\ar[r]^p&B.
      }
    \]
    Thus by comparing the homotopy exact sequence of the top and the bottom fibrations, we get that the map $\phi$ is an isomorphism in homotopy groups. Therefore the proof is done by the J.H.C. Whitehead theorem.
  \end{proof}

  Finally, we prove Theorem \ref{main 2}.

  \begin{proof}
    [Proof of Theorem \ref{main 2}]
    By the construction of the map $\mathcal{W}_L\to\Gr_L(\C^{2L})$ satisfies the commutative diagram
    \begin{equation}
      \label{W - W period}
      \xymatrix{
        \widehat{\mathcal{W}}_L\ar[r]^(.4){\rho}\ar[d]&\Gr_L(\C^{2L})\ar@{=}[d]\\
        \mathcal{W}_L\ar[r]&\Gr_L(\C^{2L})
      }
    \end{equation}
    where the map $\rho$ is as in \eqref{fibration W period}. This diagram is natural with respect to $L$, so we get a commutative diagram
    \[
      \xymatrix{
        \widehat{\mathcal{W}}\ar[r]\ar[d]&B\U(\infty)\ar@{=}[d]\\
        \mathcal{W}\ar[r]&B\U(\infty).
      }
    \]
    Theorem \ref{decomp period} implies the top map has a homotopy section, so the bottom map has a homotopy section too. Then by Lemmas \ref{deloop} and \ref{fibration trivial}, there is a homotopy equivalence
    \[
      \mathcal{W}\simeq B\U(\infty)\times\mathcal{F}.
    \]
    Since $\mathcal{U}\simeq\Z\times\mathcal{U}^0$, the proof is complete by Proposition \ref{U^0} and Lemma \ref{F}.
  \end{proof}

\end{document}